\documentclass[12pt]{article}
\usepackage{amssymb}
\usepackage{amsmath}
\usepackage{a4wide}
\usepackage{latexsym}
\usepackage{graphicx}
\newcommand {\eq}[1]{\begin{equation}\label{#1}}
\newcommand {\en} {\end{equation}}
\newcommand {\proof} {\noindent{\it Proof}. \ignorespaces}
\newcommand {\eproof}
      {\space
        {\ \vbox{\hrule\hbox{\vrule height1.3ex\hskip0.8ex\vrule}\hrule}}
        \vskip 0.3cm \par}
\newcommand {\mat}[1]{\left[\begin{array}{#1}}
\newcommand {\rix}          {\end{array}\right]}

\newtheorem{theorem}          {Theorem}
\newtheorem{lemma}         [theorem]{Lemma}
\newtheorem{definition}    [theorem]{Definition}
\newtheorem{corollary}     [theorem]{Corollary}
\newtheorem{proposition}     [theorem]{Proposition}

\newtheorem{example}        [theorem]   {Example}
\newtheorem{remark}         [theorem]   {Remark}

\newtheorem{openproblem}   [theorem] {Open Problem}
\newcommand {\diag}     {\mathop{\rm diag}\nolimits}
\newcommand {\rank}     {\mathop{\rm rank}\nolimits}

\newcommand {\size}    {\mathop{\rm size}\nolimits}

\date{\small Revision 13.10.2012}
\title{Matrices that commute with their derivative.\\
On a letter from Schur to Wielandt.
\footnotemark[2]}

\author{
Olga Holtz\footnotemark[3]
\and    Volker Mehrmann\footnotemark[3]
       \and Hans Schneider \footnotemark[4]
}

\begin{document}

\renewcommand{\thefootnote}{\fnsymbol{footnote}}

\footnotetext[3]{Institut f\"ur Mathematik, Technische Universit\"at Berlin, Stra{\ss}e des 17. Juni 136,
D-10623 Berlin, FRG.
\texttt{$\{$holtz,mehrmann$\}$@math.tu-berlin.de}.}
\footnotetext[4]{Department of Mathematics, University of
Wisconsin, Madison, WI 53706, USA. \texttt{hans@math.wisc.edu}}
\footnotetext[2]{Supported by {\it Deutsche Forschungsgemeinschaft},
through DFG project ME 790/15-2.}

 \maketitle

\begin{abstract}  We examine when a matrix
  whose elements are differentiable functions in one variable commutes with its derivative.
This problem was discussed in a letter from Issai Schur to Helmut Wielandt written in 1934, which we found in Wielandt's Nachlass. We present this letter and its translation into English.
The topic was rediscovered later and partial results were proved. However, there are  many subtle observations in Schur's letter which were not obtained in later years. Using an algebraic setting, we put these into perspective and extend them in several directions. We present in detail the relationship between several conditions mentioned in Schur's letter and we focus in particular on the characterization of matrices called  Type~1 by Schur.
We also present several examples that demonstrate Schur's observations.



\end{abstract}

\section{Introduction}

What are the conditions that force a matrix  of differentiable  functions to commute with its elementwise derivative?
This problem, discussed in a letter from I. Schur to H. Wielandt \cite{Sch34}, has been 
discussed in a large number of papers \cite{Ama54,Asc50,Asc52,BogC59,Die74,Eps63,Eru46,Gof81,Hel55,KotE82,
Kuz76,Mar65,Mar67,Par72,Pet79,Ros65,Sch52,Ter55}. However, these authors were unaware of Schur's letter and did not find some of its principal results.  A summary and a historical discussion of the problem and several extensions thereof are presented by Evard in \cite{Eva85,Eva95}, where the study of the topic is dated back to the 1940s and 1950s, but
Schur's letter shows that it already  appeared in Schur's lectures in the 1930s, if not earlier.

\if {
We do not  know which set of functions Schur had in mind,
whether Schur meant this to be the set of  analytic functions
or  meromorphic functions in one variable. He
may even have had the rational functions in one variable over the
complex numbers in mind. And we do not know which arguments Schur used to reach conclusions concerning matrices of small size at the end of his letter.  Though our arguments remain close to those of Schur, we will take an algebraic approach and discuss the results in Schur's letter in differential fields.
This is also the approach that was taken  in \cite{AdkEG93} and  in unpublished notes of Guralnick~\cite{Gur05}, where results related to ours using differential fields were discussed.
This approach has also been taken in 

} \fi

The content of the paper is as follows.
In Section~\ref{letter} we present a facsimile of Schur's letter to
Wielandt and its English translation. In Section~\ref{sec:discussion} we discuss Schur's letter and we motivate our use of differential fields.
In Section~\ref{sec:prelim} we introduce our notation and reprove Frobenius result on Wronskians. In Section~\ref{sec:ct1} we discuss the results that characterize
the matrices of Type~1 in Schur's letter and in our main Section~\ref{sec:triandia} we discuss
the role  played by diagonalizability and triangularizability of the matrix  in the commutativity of the matrix and its derivative.
We also present several illustrative examples in Section~\ref{sec:exs}
and we state an open problem in Section~\ref{conclusion}.

\section{A letter from Schur to Wielandt}\label{letter}
Our paper deals with the following letter from Issai Schur to his PhD
student Helmut Wielandt. See the facsimile below.
%
\begin{figure}
\scalebox{.8}{\includegraphics{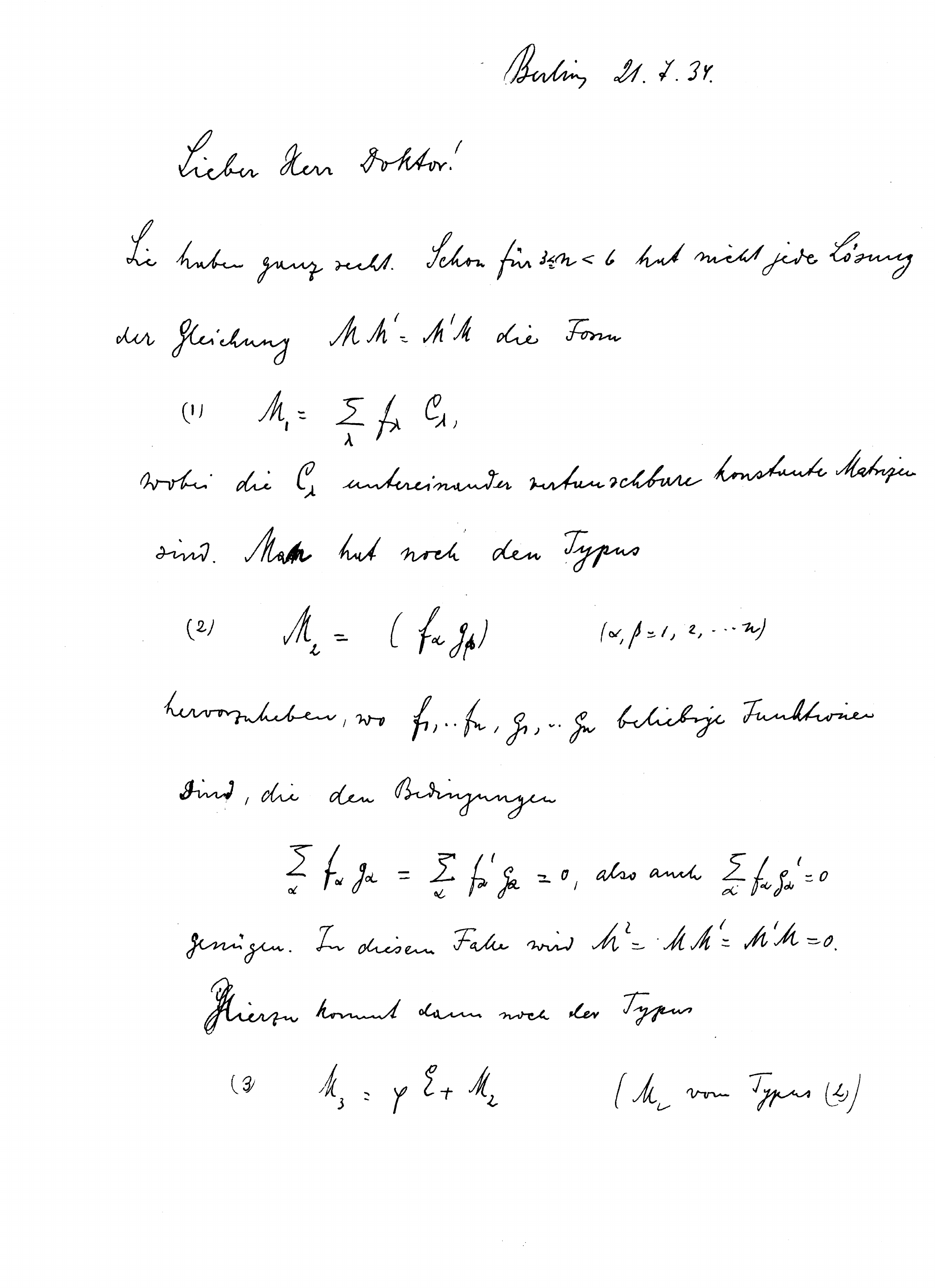}}
 \vspace{0.5cm}
	\label{FigOne}
\end{figure}
\begin{figure}
\scalebox{.8}{\includegraphics{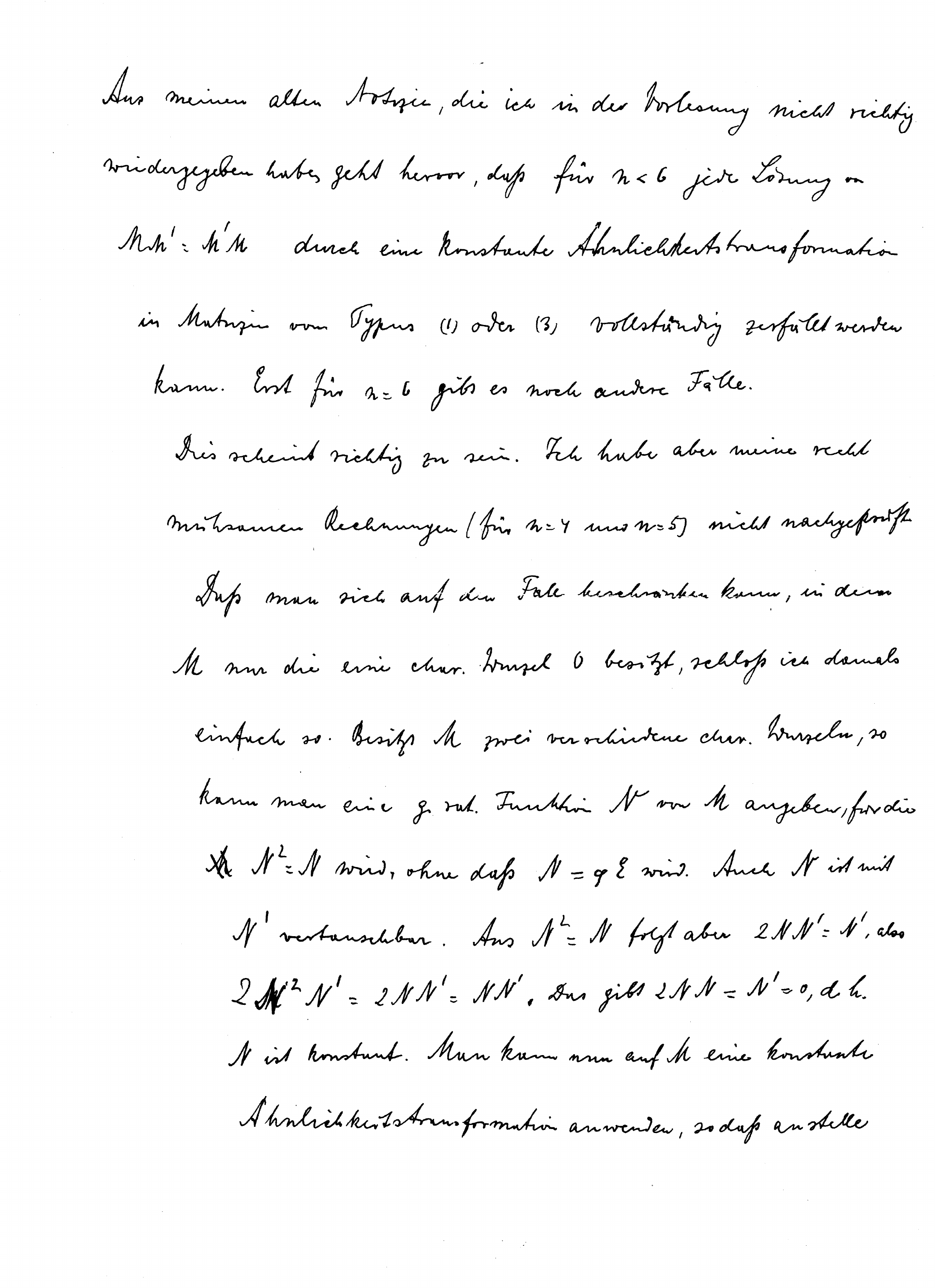}}
 \vspace{0.5cm}
	\label{FigTwo}
\end{figure}
\begin{figure}
\scalebox{.8}{\includegraphics{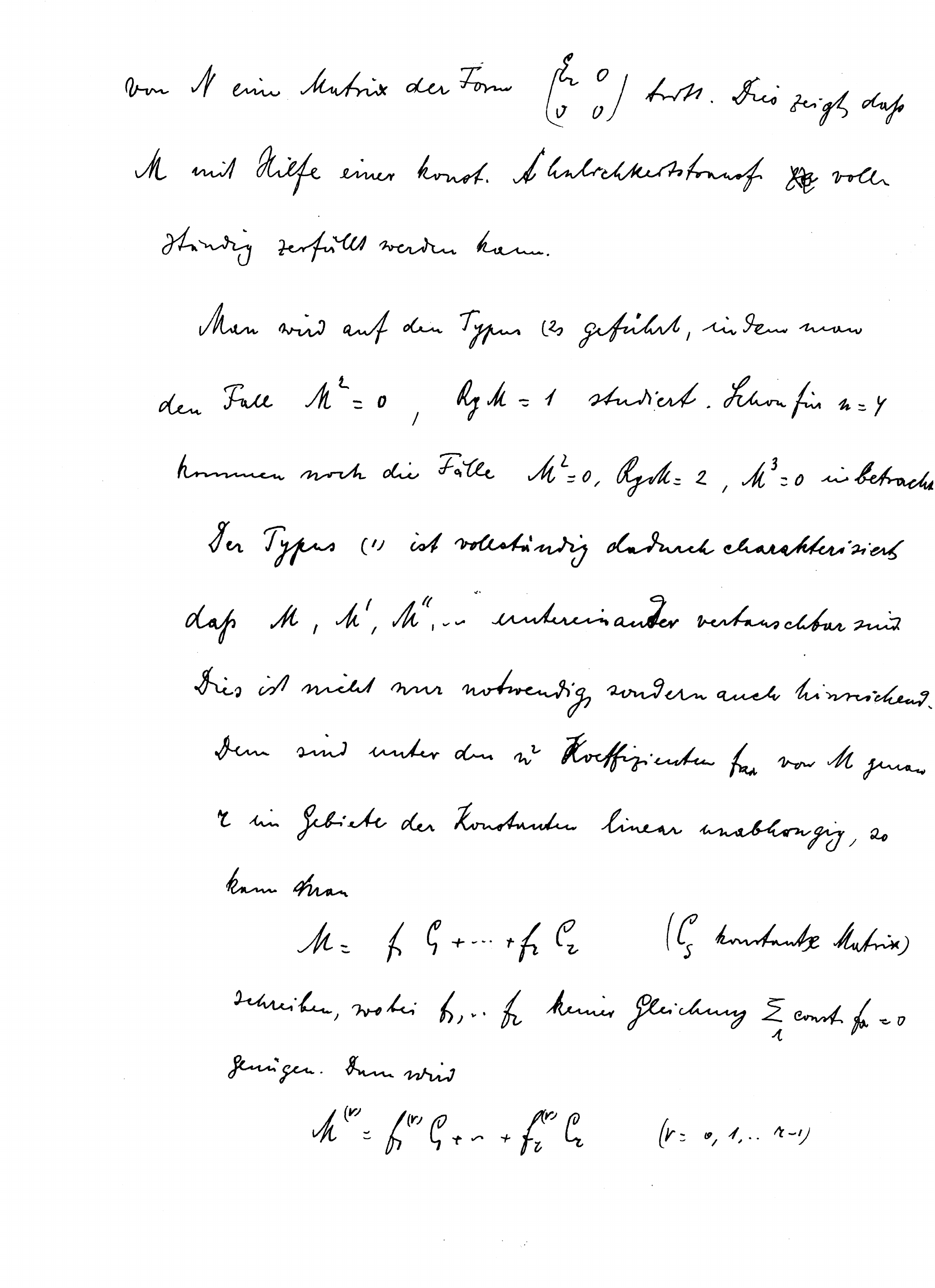}}
 \vspace{0.5cm}
	\label{FigThree}
\end{figure}
\begin{figure}
\scalebox{.8}{\includegraphics{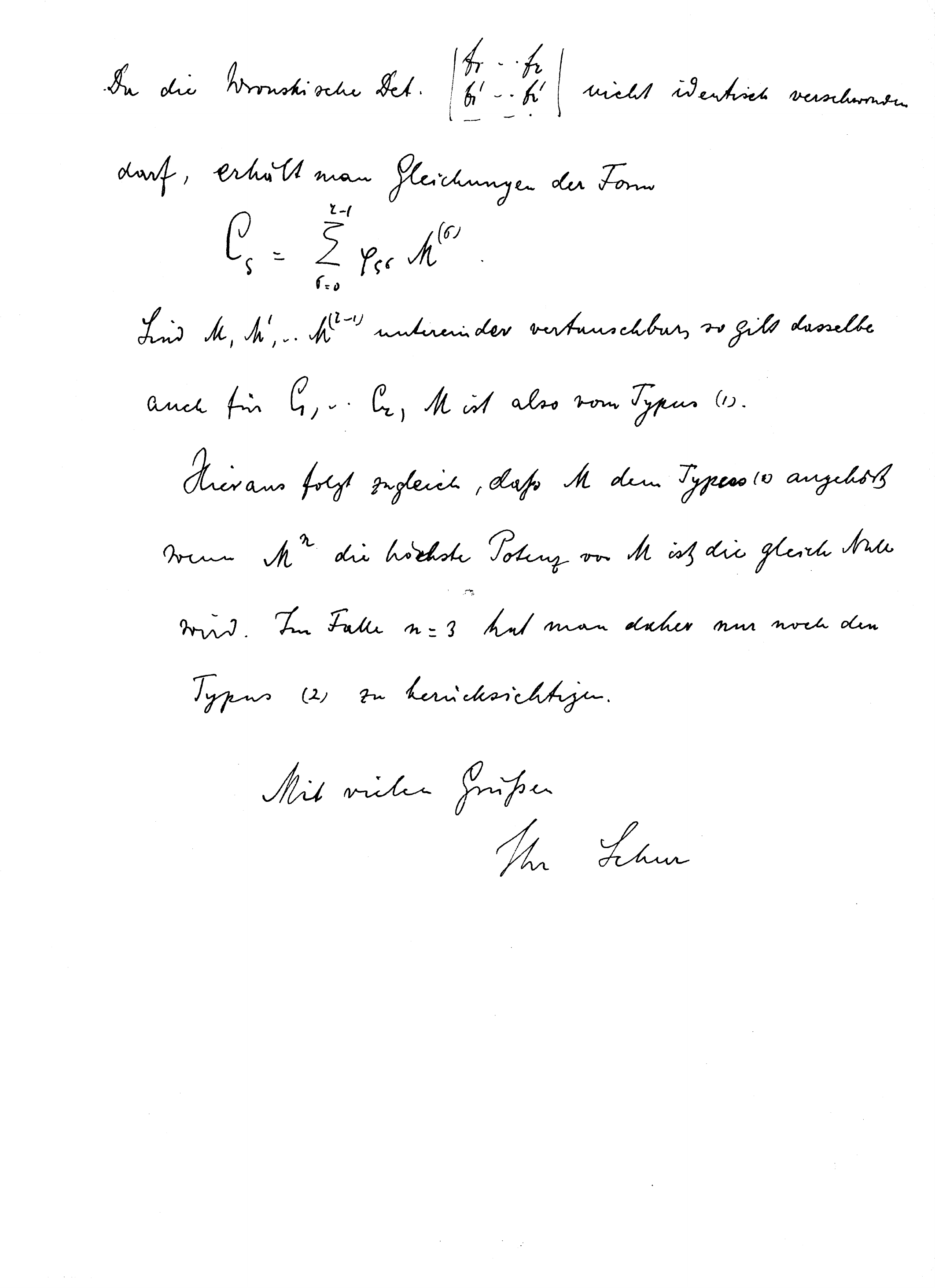}}
 \vspace{0.5cm}
	\label{FigFour}
\end{figure}
%
Translated into English, the letter reads as follows:

{\it
Lieber Herr Doktor!\hfill Berlin, 21.7.34

You are perfectly right. Already for $3\leq n<6$
not every solution of the
equation $M M'=M'M$ has the form
\eq{1} M_1 =\sum_{\lambda} f_\lambda C_\lambda, \en
where the $C_\lambda$ are pairwise commuting constant matrices. One
must also consider the type
\eq{2} M_2=(f_\alpha g_\beta), \quad (\alpha,\beta=1,\ldots n), \en
where $f_1,\ldots f_n$, $ g_1,\ldots,g_n$ are arbitrary functions
that satisfy the conditions
$$ \sum_\alpha f_\alpha g_\alpha= \sum_\alpha f'_\alpha g_\alpha=0$$
and therefore also
$$ \sum_\alpha f_\alpha g'_\alpha=0.$$
In this case we obtain
$$ M^2 =M M'=M'M=0.$$
In addition we have the type
\eq{3} M_3=\phi E+M_2, \en
with $M_2$ of type (\ref{2}). \footnote{Note that $E$ here denotes the identity matrix.}
 From my old notes, which  I did not present correctly in my lectures,
it can be deduced that for $n<6$ every solution of $M M'=M'M$
can be completely decomposed by means of constant similarity
transformations into matrices of type (\ref{1}) and (\ref{3}).
Only from $n=6$ on there are also other cases.
This seems to be correct.
But I have not checked my rather laborious
computations (for $n=4$ and $n=5$).

I concluded in the following simple manner that one can restrict oneself to
the case where $M$ has only one characteristic root (namely $0$):  If $M$
has two different characteristic roots, then one can determine
a  rational  function  $N$ of $M$ for which $N^2=N$ but
not $N=\phi E$. Also $N$ commutes with $N'$. It follows from $N^2=N$
that $2NN'=N'$, thus $2N^2N'= 2NN'=NN'$. This yields
$2 NN'=N'=0$, i.e., $N$ is constant.

Now one can apply a constant similarity transformation to $M$ so that
instead of $N$ one achieves a matrix of the form
$$\mat{cc} E & 0 \\ 0 & 0 \rix.$$
This shows that $M$ can be decomposed  completely
by means of a constant
similarity transformation.

One is led to type (\ref{2}) by studying the case $M^2=0$, $\rank M=1$.
Already for $n=4$ also the cases $M^2=0$, $\rank M=2$, $M^3=0$ need to be
considered.

Type (\ref{1}) is completely characterized by the property
that $M,M',M'', \ldots $ are pairwise
commuting. This is not only necessary but also sufficient.
For, if among the $n^2$ coefficients $f_{\alpha \beta}$ of $M$ exactly
$r$ are linearly independent over the domain of constants,
then one can write
$$ M=f_1 C_1+\cdots +f_rC_r,$$
 ($C_s$  a constant matrix), where $f_1, \ldots, f_r$ satisfy no
 equation $\sum_{\alpha} \mbox{\rm const}\ f_\alpha=0$.
Then
$$ M^{(\nu)}= f_1^{(\nu)} C_1+ \cdots +f_r^{(\nu)} C_r,
\quad (\nu=1,\ldots,r-1).$$
Since the Wronskian determinant
$$\left |\mat{ccc} f_1 &\ldots& f_r \\
f_1' &\ldots& f_r'\\
 & \vdots & \rix\right |
 $$
 cannot vanish identically, one obtains equations of the form
$$ C_s=\sum_{\sigma=0}^{r-1} \phi_{s \sigma} M^{(\sigma)}.$$
If $M,M',M'', \ldots, M^{(r-1)}$ are pairwise
commuting, then the same is true also for
$C_1,\ldots C_r$ and thus $M$ is  of type (\ref{1}). This implies
furthermore that $M$ belongs to type (\ref{1}) if
$M^n$ is the highest\footnote{We think that Schur means \emph{lowest} here.} power of $M$ that equals $0$.
In the case $n=3$ one therefore only needs to consider type (\ref{2}).

\hfill With best regards

\hfill Yours, Schur
}

\bigskip

\section{Discussion of Schur's letter}\label{sec:discussion}

This letter was found in Helmut Wielandt's mathematical Nachlass when it was collected by  Heinrich Wefelscheid and Hans Schneider not long after Wielandt's death in 2001.  We may therefore safely assume that Schur's recent student Wielandt is the "Herr Doktor" to whom the letter is addressed. Schur's letter begins with a  reference to a previous  remark of Wielandt's which corrected an incorrect assertion by Schur. We can only guess at this sequence of events, but perhaps a clue is provided by Schur's reference to his notes  which he did not present correctly in his lectures. Could Wielandt have been in the audience and did he subsequently point out the error? And what was this error?  Very probably it was that every matrix of functions that commutes with its derivative is given by (1) (matrices called Type 1), for Schur now denies this and displays another type of matrix commuting with its derivative (called Type 2). He recalls that in his notes he claimed that for matrices of size 5 or less every such matrix is of Type 1, 2 or 3, where Type 3 is obtained from Type 2 by adding a scalar function times the identity. This is not correct because there is also the direct sum of a size 2 matrix of Type 1 and a size 3 matrix of Type 2, we prove this below.

We do not know why Schur was interested in the topic of matrices of functions that commute with their derivative, but
it is probably a safe guess that this question came up in the context of solving differential equations, at least this is
the motivation in many of the subsequent papers on this topic.

As one of the main results of his letter, Schur shows that an idempotent that commutes with its derivative is a constant matrix  and, without further explanation,  concludes that one can restrict oneself to matrices with a single eigenvalue. The latter  observation raises several questions. First,  Schur does not say which functions he has in mind. Second, his argument follows from a standard decomposition of a matrix by a similarity into a direct sum of matrices {\em provided that} the eigenvalues of the matrix are functions of the type considered.  But this is not true in general, for example the eigenvalues of a matrix of rational functions are algebraic functions.
We wonder whether Schur was aware of this difficulty and we shall return to it at the end of this section.

Then Schur shows a  matrix of size $n$ is of Type 1 if and only if it and its first $n-1$ derivatives are pairwise commutative. His proof is based on a result of Frobenius \cite{Fro1874} that a set of functions is linear independent over the constants if and only if their Wronskian determinant is nonzero. Frobenius, like Schur,  does  not explain what functions he has in mind. In fact,  Peano \cite{Pea1889} shows that there exist real differentiable functions that are linearly independent over the reals whose Wronskian is $0$. This is followed by Bocher \cite{Boc00} who shows that Frobenius' result  holds for analytic functions and investigates necessary and sufficient conditions in \cite{Boc01}. A good discussion of this topic can be found in \cite{BosD10}.

We conclude this section by explaining how our exposition has been influenced by some of the observation above. As we do not know what functions Schur and Frobenius had in mind, we follow \cite{AdkEG93} and some unpublished notes  of Guralnick \cite{Gur05} and set Schur's results and ours in terms of differential fields (which include the field of rational functions and the quotient field of analytic functions over the real or complex numbers).  Since we do not know how Schur concludes that it is enough to consider matrices with a single eigenvalue, we derive our results from  standard matrix decomposition (our Lemma \ref{decomp} below) which does not assume that all eigenvalues lie in the differential field under consideration.

 \section{Notation and preliminaries}\label{sec:prelim}

A {\em differential field} $\mathbb F$ is an (algebraic) field together
with an additional operation (the derivative), denoted by ${}'$ that satisfies $(a+b)' = a'
+ b'$ and $(ab)' = ab' + a'b$ for $a,b \in \mathbb F$. An element $a
\in \mathbb F$ is called a {\em constant} if $a' = 0$. It is easily
shown that the set of constants forms a subfield $\mathbb K$ of $\mathbb F$
with $1 \in \mathbb K$. Examples are provided by the rational functions over the real or complex numbers
and the meromorphic functions over the complex numbers.

In what follows we consider a (differential) field $\mathbb F$
and  matrices $M=[m_{i,j}]\in \mathbb F^{n,n}$.
The main condition that we want to analyze is when $M\in \mathbb F^{n,n}$
commutes with its derivative,
\begin{equation}\label{c1}
MM'=M'M.
\end{equation}
As  $M\in \mathbb F^{n,n}$, it has a minimal and a characteristic polynomial,  and $M$ is called {\em nonderogatory} if the characteristic polynomial is equal to the minimal  polynomial, otherwise it is called {\em derogatory}. See \cite{HorJ85}. 

In Schur's letter the following
 three types of matrices are considered.
\begin{definition}
Let $M\in \mathbb F^{n,n}$. Then $M$ is said to be of
\begin{itemize}
\item {\em Type 1\/} if
\[
M =\sum_{j=1}^{k} f_j C_j,
\]
where $f_j\in \mathbb F$, and $C_j\in \mathbb K^{n,n}$, for $j=1,\ldots,k$,
and the $C_j$ are pairwise commuting;
\item {\em Type 2\/} if
\[
M=f g^T,\
\]
with $f,g\in \mathbb F^{n}$, satisfying $f^Tg=f^Tg'=0$;
\item {\em Type 3\/} if
\[
M=hI+\widetilde M,
\]
with $h\in \mathbb F$ and $\widetilde M$ is of Type~2.
\end{itemize}
\end{definition}
%
Schur's letter also mentions the condition that all derivatives of $M$ commute, i.e.,
\begin{equation}\label{c6}
M^{(i)}M^{(j)}=M^{(j)}M^{(i)}\ \mbox{\rm for all nonnegative integers}\ i,j.
\end{equation}
%
%
%
%
%
%
To characterize the relationship between all these properties, we first recall
several results from Schur's letter and from classical algebra.
\begin{lemma} \label{insight} Let $\mathbb F$ be a differential field with field of constants
 $\mathbb K$. Let $N$ be an idempotent matrix in
$\mathbb {F}^{n,n}$ that commutes with $N'$. Then $N \in \mathbb K^{n,n}$.
\end{lemma}
\proof (See Schur's letter.) It follows from $N^2 = N$ that $2NN' =
N'$. Thus $2NN' =  2N^2N' = NN'$ and this implies that $0 = 2NN' = N'$.
\eproof
Another important tool in our analysis will be the following result which in its original form
is due to Frobenius~\cite{Fro1874}, see Section~\ref{sec:discussion}.
We phrase and prove the result in the context of differential fields.
\begin{theorem}\label{frob}
Consider a differential field $\mathbb F$ with field of constants $\mathbb K$. Then
$y_1,\ldots,y_r\in \mathbb F$ are linearly dependent over $\mathbb K$ if and only if the columns of the \emph{Wronski matrix}
\[
Y=\left[ \begin{array}{cccc} y_1 & y_2 & \dots & y_r \\
                         y_1' & y_2' &\dots & y_r'\\
                         \vdots & \vdots & \ddots & \vdots\\
                         y_1^{(r-1)} & y_2^{(r-1)} &\dots & y_r^{(r-1)}
\end{array} \right ],
\]
are linearly dependent over  $\mathbb F$.
\end{theorem}
\proof
We proceed by induction over $r$. The case $r=1$ is trivial.

Consider the Wronski matrix $Y$ and the lower triangular matrix
  \[
Z=\left[ \begin{array}{cccc} z & 0 & \dots & 0 \\
                         c_{2,1}z' & z &\dots & 0\\
                         \vdots & \vdots & \ddots & \vdots\\
                         c_{n,1}z^{(n-1)} & c_{n,2}z^{(n-2)} &\dots & z
\end{array} \right ],
\]
with $c_{i,j}$ appropriate binomial coefficients such that
\[
ZY=\left[ \begin{array}{cccc} z  y_1&   z  y_2& \dots & z  y_r \\
                         (z  y_1)' & (z  y_2)' &\dots & (z  y_r)'\\
                         \vdots & \vdots & \ddots & \vdots\\
                         (z  y_1)^{(n-1)} & (z  y_2)^{(n-1)} &\dots & (z  y_r)^{(n-1)}
\end{array} \right ].
\]
Since $\mathbb F$ is a differential field, we can choose $z=y_1^{-1}$ and obtain
that
\[
ZY=\left[ \begin{array}{cccc} 1&   y_1^{-1}  y_2& \dots & y_1^{-1}  y_r \\
                         0 & (y_1^{-1} y_2)' &\dots & ( y_1^{-1} y_r)'\\
                         \vdots & \vdots & \ddots & \vdots\\
                         0 & (y_1^{-1}  y_2)^{(n-1)} &\dots & ( y_1^{-1} y_r)^{(n-1)}
\end{array} \right ].
\]
It follows that the columns of $Y$ are linearly dependent over $\mathbb F$ if and only if the columns
of
\[
\left[ \begin{array}{ccc}
                          (y_1^{-1} y_2)' &\dots & ( y_1^{-1} y_r)'\\
                         \vdots  & \ddots & \vdots\\
                          (y_1^{-1}  y_2)^{(n-1)} &\dots & ( y_1^{-1} y_r)^{(n-1)}
\end{array} \right ]
\]
are linearly dependent over $\mathbb F$, which, by induction, holds
 if and only if  $(y_1^{-1} y_2)' ,\dots , ( y_1^{-1} y_r)'$
are linearly dependent over $\mathbb K$, i.e., there exist coefficients
$b_2,\ldots,b_r\in \mathbb K$, not all $0$,
such that
\[
b_2 (y_1^{-1} y_2)'+\cdots +b_r( y_1^{-1} y_r)'=0.
\]
Integrating this identity, we obtain
\[
b_2 (y_1^{-1} y_2)+\cdots +b_r( y_1^{-1} y_r)=-b_1
\]
for some integration constant $b_1\in \mathbb K$, or equivalently
\[
b_1 y_1+\cdots + b_r y_r=0.
\]
{}\vskip -1cm\hskip 12cm  \eproof

Theorem~\ref{frob} implies in particular that the columns of the Wronski matrix $Y$
are linearly independent over $\mathbb F$ if and only if they are linearly independent over $\mathbb K$.

\begin{remark}{\rm
Theorem~\ref{frob} is discussed from a formal algebraic point of view, which however includes the cases of complex analytic functions and rational functions over a field, since these are contained in differential fields. Necessary and sufficient conditions  for Theorem~\ref{frob}  to hold for other functions were proved in \cite{Boc01} and discussed in many places, see, e.g., \cite{BosD10,Mei61} and \cite[Ch. XVIII]{Mui33}.
 }
\end{remark}

\section{Characterization of matrices of Type~1} \label{sec:ct1}
In this section we discuss  relationships among the various properties introduced in Schur's letter and in the previous section. This will give, in particular, a characterization of matrices of Type~1.

In his letter, Schur proves the following result.

\begin{theorem}\label{th:allcommute}
Let $\mathbb F$ be a differential field. Then $M\in {\mathbb F}^{n,n}$
is of Type~1 if and only if it satisfies condition (\ref{c6}), i.e., $M^{(i)}M^{(j)}=M^{(j)}M^{(i)}$ for all nonnegative integers $i,j$.
\end{theorem}
\proof (See Schur's letter.)
If $M$ is of Type~1, then $M =\sum_{j=1}^{k} f_j C_j$ and the $C_j\in \mathbb K^{n,n}$ are pairwise commuting, which immediately implies (\ref{c6}).  For the converse,
Schur makes use of Theorem~\ref{frob}, since if among the $n^2$ coefficients
$m_{i,j}$ exactly $r$ are linearly independent over $\mathbb K$, then
\[
M=f_1C_1+\cdots+f_rC_r,
\]
with coefficients $C_i\in \mathbb K^{n,n}$, where
$f_1,\ldots,f_r$ are linearly independent over $\mathbb K$.
Then
\[
M^{(i)}=f_1^{(i)}C_1+\cdots +f_r^{(i)}C_r,\qquad i=1,\ldots,r-1.
\]
By Theorem~\ref{frob}, the columns of the associated Wronski matrix
are linearly independent, and hence each of the $C_i$ can be expressed
as
\[
C_i=\sum_{j=0}^{r-1} g_{i,j} M^{(j)}.
\]
Thus, if condition~(\ref{c6}) holds, then the $C_i$, $i=1,\ldots,r$,
are pairwise commuting and thus $M$ is of Type~1.
\eproof
Using this result we immediately have the following Theorem.
\begin{theorem}\label{th:nondero}

Let $\mathbb F$ be a differential field with field of constants
$\mathbb K$. If $M\in \mathbb F^{n,n}$ is nonderogatory  and  $MM' = M'M$,
then $M$ is of Type~1.
\end{theorem}
\proof
If $M$ is nonderogatory then all matrices that commute with $M$ have the form $p(M)$, where $p$ is a polynomial with coefficients in $\mathbb F$, see \cite{DraDG51,HorJ85}.  Thus $MM' = M'M$ implies that $M'$ is a polynomial in $M$.
But then every derivative $M^{(j)}$ is a polynomial in $M$ as well and
thus (\ref{c6}) holds which by Theorem~\ref{th:allcommute} implies that $M$ is of Type~1.
\eproof
The following example from \cite{Asc52,Eva85} of a Type~2 matrix shows that one cannot easily drop the
condition that the matrix is nonderogatory.

\begin{example}\label{eva3x3}
{\rm   Let
\[ f =\mat{c}   1\\   t \\ t^2\rix,\qquad g =\mat{c}  t^2\\ -2 t\\   1
\rix,
\]
then $f^T g=0$ and  $f^T g'=0$, hence
\begin{equation}\label{evaM}
M_a := g f^T= \mat{ccc} t^2 &   t^3&   t^4\\
-2 t& -2 t^2 & -2 t^3\\
1&     t &  t^2\rix,
\end{equation}
is of Type~2. Since $M_a$ is nilpotent with $M_a^2=0$ but $M_a\neq 0$ and the rank is $1$,
it is derogatory. One has
\[
M_a'=
 \mat{ccc}    2t & 3t^2&  4t^3\\
-2  & -4t& -6t^2\\
0 &      1 &    2 t
\rix,\qquad
M_a''= \mat{ccc}   2 & 6t & 12 t^2 \\
 0 & -4 & -12 t \\
0 & 0 & 2\rix,
\]
and thus $ M_a M_a'=M_a'  M_a=0$. By the product rule
it immediately follows that $M_a M_a''=M_a''M_a$, but
\[
M_a' M_a''= \mat{ccc}
  4t&     0& -4t^3\\
    -4  & 4t &12t^2\\
     0  &   -4 &  -8t
\rix \neq M_a'M_a''=\mat{ccc}
-8t &-6t^2 & -4t^3\\
   8  &  4t  &    0\\
   0  &   2 &  4t\rix.
\]
Therefore, it follows from Theorem~\ref{th:allcommute} that $M_a$ is not of Type~1.

For any dimension $n\geq 3$, one can construct an example of Type~2 by choosing $f \in {\mathbb F}^n$, setting $F = [f,f']$ and then choosing $g$ in the nullspace of $F^T$.
Then  $fg^T$ is of Type~2.
}
\end{example}
Actually every nilpotent matrix function $M$ of rank one
satisfying $MM' = M'M$ is of the form $M=fg^T$ and hence of Type~2.
This follows immediately because if $M=fg^T$ and $M^2=0$ then $g^Tf=0$
and hence $g^Tf'+(g^T)'f=0$. Then it follows from $MM' = M'M$
that $fg^T(f(g^T)' + f'g^T)  = (g^Tf')fg^T= (f(g^T)' + f'g^T)fg^T  = (g^T)'f fg^T$
which implies that  $g^Tf' = f^Tg'$ and hence $g^Tf' = f^Tg' = 0$.

\section{Triangularizability and Diagonalizability}\label{sec:triandia}

In his letter Schur claims that it is sufficient to consider the case
that $M\in \mathbb F^{n,n}$ is triangular with only one eigenvalue. This  follows from his argument
in case the matrix has its eigenvalues in $\mathbb F$, which could be guaranteed by assuming that this
matrix is $\mathbb F$-diagonalizable or even $\mathbb F$-triangularizable. Clearly a sufficient condition
for this to hold is that $\mathbb F$ is algebraically closed, because then for every matrix in $\mathbb F^{n,n}$ the characteristic polynomial splits into linear factors.
%
\if{
To do this, an essential property discussed in Schur's letter is whether there exists a similarity
transformation to triangular or diagonal form with nonsingular matrices in $\mathbb F^{n,n}$ or
$\mathbb K^{n,n}$.
}\fi
%
\begin{definition}\label {triangularizable}
Let $\mathbb F$ be a differential field and let $\mathbb H$ be a subfield of $\mathbb F$.
Then $M\in \mathbb F^{n,n}$ is called
{\em $\mathbb H$-triangularizable (diagonalizable)} if there exists a nonsingular $T\in \mathbb H^{n,n}$ such that $T^{-1} M T$ is upper triangular (diagonal).
\end{definition}
Using Lemma~\ref{insight}, we can obtain the following result for
matrices $M\in \mathbb F^{n,n}$ that commute with their derivative $M'$, which is most likely well known
but we could not find a reference.
\begin{lemma}\label{decomp}
Let $\mathbb F$ be a differential field with field of constants $\mathbb K$, and suppose that $M \in {\mathbb F}^{n,n}$ satisfies $MM' = M'M$. Then there
exists an invertible matrix $T\in \mathbb K^{n,n}$ such that
\eq{dirsum}
T^{-1} M T= \diag (M_1,\ldots,M_k),
\en
where the minimal polynomial of each $M_i$ is a power of a
polynomial that is irreducible over $\mathbb F$.
\end{lemma}
\proof
 Let the minimal polynomial of $M$ be $\mu(\lambda) = \mu_1(\lambda)\cdots \mu_k(\lambda)$,
where the $\mu_i(\lambda)$ are powers of pairwise distinct polynomials that are
irreducible over $\mathbb F$. Set
\[
p_i(\lambda) = \mu(\lambda)/\mu_i(\lambda),\qquad i = 1,\ldots,k.
\]
Since the polynomials $p_i(\lambda)$ have no common factor, there exist polynomials
$q_i(\lambda)$, $i = 1,\ldots,k$, such that the polynomials $\epsilon_i(\lambda) = p_i(\lambda)q_i(\lambda)$, $i = 1,\ldots,k$, satisfy
\begin{equation}\label{epscomp}
 \epsilon_1(\lambda) + \cdots + \epsilon_k(\lambda) = 1.
\end{equation}
Setting $E_i = \epsilon_i(M)$, $i = 1, \ldots, k$ and using the fact that
$\mu(M) =0$ yields that
\begin{eqnarray}\label{comp}
E_1 + \cdots + E_k &=& I, \\
\label{ortho}
E_iE_j &=& 0,\qquad 
\; \, i,j = 1,\ldots,k,\quad  i\neq j,\\
\label{idem}
E_i^2 &=& E_i,\qquad i = 1,\ldots,k.
\end{eqnarray}
The last identity follows directly from (\ref{comp}) and (\ref{ortho}).
Since the $E_i$ are polynomials in $M$ and $MM' = M'M$, it follows that  the $E_i$ commute with $E'_i$, $i = 1,\ldots k$.  Hence, by Lemma~\ref{insight}, $E_i \in \mathbb K^{n,n}$, $i = 1,\ldots,k$. By (\ref{comp}), (\ref{ortho}), and (\ref{idem}), $\mathbb K^n$ is a direct sum
of the ranges of the $E_i$ and we obtain  that, for some nonsingular  $T\in \mathbb K^{n,n}$,
\[
  \widetilde{E_i}:= T^{-1}E_iT = \diag(0, I_i,0),\qquad   i = 1,\ldots,k,
\]
where the $I_i$ are identity matrices of the size equal to the dimension to the
range of $E_i$. This is a consequence of the fact that $E_i$ is diagonalizable with eigenvalues $0$ and $1$.
Since each $E_i$ commutes with $M$,
we obtain that
\begin{eqnarray*}
 \widetilde{M_i}&:=& T^{-1}E_iMT \\
 &=&  T^{-1}E_iME_iT\\
 &=& \diag(0,I_i,0) T^{-1}MT \diag(0,I_i,0)\\
 &=&\diag(0,M_i,0) ,\qquad  i = 1, \ldots, k.
\end{eqnarray*}
Now observe that
\[
\widetilde{E_i}\mu_i(\widetilde{M}_i) \widetilde{E_i}= T^{-1} \epsilon_i(M)\mu_i(M)  \epsilon_i(M)T = 0,
\]
since
$\epsilon_i(\lambda)\mu_i(\lambda )= \mu(\lambda)q_i(\lambda)$. Hence $\mu_i(M_i)=0$ as well.
We assert that  $\mu_i(\lambda)$ is the minimal polynomial of $M_i$, for if
 $r(M_i) = 0$ for a proper factor $r(\lambda)$ of $m_i(\lambda)$ then $r(M)\Pi_{j \neq i} \mu_j(M) = 0$,
 contrary to the assumption that $ \mu(\lambda)$ is the minimal polynomial of $M$.
 \eproof
Lemma~\ref{decomp} has the following corollary, which has been proved in a different way
in \cite{AdkEG93} and \cite{Gur05}.
\begin{corollary}\label{th:diagonal}
Let $\mathbb F$ be a differential field with field of constants
$\mathbb K$. If $M\in \mathbb F^{n,n}$ satisfies $MM' = M'M$ and is $\mathbb F$-diagonalizable, then $M$ is $\mathbb K$-diagonalizable.
\end{corollary}
\proof
In this case, the minimal polynomial of $M$ is a product of distinct linear factors and hence, the minimal polynomial of each $M_i$ occurring in the proof of Lemma~\ref{decomp} is linear. Therefore, each $M_i$ is a scalar matrix.
\eproof
We also have the following Corollary.
\begin{corollary}\label{cor:diag_type1}
Let $\mathbb F$ be a differential field with field of constants
$\mathbb K$. If $M\in \mathbb F^{n,n}$ satisfies $MM' = M'M$ and is $\mathbb F$-diagonalizable, then $M$ is of Type~1.
\end{corollary}
\proof
By Corollary~\ref{th:diagonal}, $M=T^{-1} \diag (m_1,\ldots,m_n)T$ with $m_i\in \mathbb F$
and nonsingular $T\in \mathbb K^{n,n}$. Hence
\[
M=\sum_{i=1}^n m_{i}T^{-1} E_{i,i} T
\]
where $E_{i,i}$ is a matrix that has a $1$ in position $(i,i)$ and zeros everywhere else. Since all the matrices $E_{i,i}$ commute,
$M$ is of Type~1.
\eproof
%
%

\begin{remark}\label{refrem} {\rm Any $M\in \mathbb F^{n,n}$ that is of rank one, satisfies $MM' = M'M$ and is not nilpotent, is of Type~1, since in this case $M$ is $\mathbb F$-diagonalizable. This follows by Corollary~\ref{cor:diag_type1}, since
the minimal polynomial has the from $(\lambda-c)\lambda$ for some $c\in \mathbb F$. This means in particular for a rank one matrix $M\in \mathbb F^{n,n}$ to be of Type~2 and not of Type~1 it has to be nilpotent.}
\end{remark}
For matrices that are just triangularizable the situation is more subtle. We have the following theorem.
\begin{theorem} \label{th:triangularize} Let $\mathbb F$ be a differential field with an algebraically closed field of constants $\mathbb K$. If $M \in \mathbb F^{n,n}$ is Type~1, then $M$ is $\mathbb K$-triangularizable.
\end{theorem}
\proof
Any finite set of pairwise
commutative matrices with elements in an algebraically closed field may be simultaneously triangularized, see e.g.,
\cite[Theorem 1.1.5]{RadR00}. Under this assumption on $\mathbb K$, if $M$ is Type~1, then it follows that  the matrices $C_i \in \mathbb K^{n,n}$  in the representation of $M$ are simultaneously triangularizable by a matrix
$T \in \mathbb K^{n,n}$. Hence $T$ also triangularizes $M$.
\eproof
Theorem~\ref{th:triangularize} implies that  Type~1 matrices
have $n$ eigenvalues  in $\mathbb F$ if $\mathbb K$ is algebraically closed and it further immediately
leads to a Corollary of Theorem~\ref{th:nondero}.
\begin{corollary}\label{cor:triangularize}
Let $\mathbb F$ be a differential field with field of constants
$\mathbb K$. If $M\in \mathbb F^{n,n}$ is nonderogatory, satisfies $MM' = M'M$ and
if $\mathbb K$ is algebraically closed, then $M$ is $\mathbb K$-triangularizable.
\end{corollary}
\proof
By Theorem~\ref{th:nondero} it follows that $M$ is Type~1 and thus the assertion follows
from Theorem~\ref{th:triangularize}.
\eproof

\section{Matrices of small size and examples}\label{sec:exs}

Example~\ref{eva3x3} again shows that it is difficult to drop
some of the assumptions, since this matrix is derogatory, not of Type~1, and not $\mathbb K$-triangularizable.

One might be tempted to conjecture that any $M\in \mathbb F^{n,n}$ that is $\mathbb K$-triangularizable and satisfies (\ref{c1}) is of Type~1 but this is so only for small dimensions and
is no longer true for large enough $n$, as we will demonstrate below.
Consider small dimensions  first. 
\begin{proposition}\label{2x2}
Consider a differential field $\mathbb F$ of functions with field of constants
$\mathbb K$.
Let $M=[m_{i,j}]\in {\mathbb F}^{2,2}$ be upper triangular and satisfy $M\, M'=M'\, M$. Then 
$M$ is of Type~\ref{1}.
\end{proposition}
\proof
Since $M\, M'=M'\, M$  we obtain
\[
 m_{1,2}(m_{1,1}'-m_{2,2}')-m_{1,2}'(m_{1,1}-m_{2,2})=0,
 \]
which implies that  $m_{1,2}=0$ or $m_{1,1}-m_{2,2}=0$ or
 both are nonzero and ${ d \over dt}({m_{1,1}-m_{2,2}\over m_{1,2}})=0$, i.e.,
 $c m_{1,2}+ (m_{1,1}-m_{2,2})=0$ for some nonzero constant $c$.

If $m_{1,1}=m_{2,2}$ or $m_{1,2}=0$,  then $M$, being triangular, is obviously of Type~1.
Otherwise
\[
M=m_{1,1} I + m_{1,2}\mat{cc} 0 & 1 \\ 0 & c \rix.
\]
and hence again of Type~1 as claimed.
\eproof
Proposition~\ref{2x2} implies  that $2\times2$ $\mathbb K$-triangularizable matrices
satisfying (\ref{c1}) are of Type~1. 

\begin{proposition}\label{2x2type}
Consider a differential field $\mathbb F$ with an algebraically closed field of constants $\mathbb K$.
Let $M=[m_{i,j}]\in {\mathbb F}^{2,2}$ satisfy $M\, M'=M'\, M$.
Then $M$ is of Type~1.
\end{proposition}
\proof If $M$ is $\mathbb F$-diagonalizable, then the result follows by Corollary~\ref{cor:diag_type1}. If $M$ is not $\mathbb F$-diagonalizable, then it is nonderogatory and the result follows by Corollary~\ref{cor:triangularize}.
\eproof

\begin{example}{\rm
In the $2\times2$ case, any Type~2 or Type~3 matrix is also of Type~1 but not every
Type~1 matrix is Type~3.

Let $M=\phi I_2 +f g^T$  with
\[\phi\in \mathbb F,\quad f=\mat{c} f_1 \\ f_2\rix ,\quad g=\mat{c} g_1 \\ g_2\rix\in {\mathbb F}^2\]
 be of Type~3, i.e.,
$f^Tg={f'}^Tg= f^T{g'}=0$.

If $f_2=0$, then $M$ is upper triangular and hence by Proposition~\ref{2x2}, $M$ is of Type~\ref{1}.
If $f_2\neq 0$, then with
\[
T=\mat{cc} 1 & -f_1/f_2 \\ 0 &  1 \rix.
\]
we have
\[
T MT^{-1}= \phi I_2 +\mat{cc} 0 & 0 \\ f_2 g_1 & 0 \rix=
\phi I_2 +f_2g_1\mat{cc} 0 & 0 \\ 1 & 0 \rix,
\]
since  $ f_1g_1+f_2 g_2=0$, and hence $M$ is  of Type~1.

However,  if we consider
\[
M=\phi I_2+ f\mat{cc} 0 & c \\ 0 & d\rix
\]
with $\phi,f$ nonzero functions and $c,d$ nonzero  constants,
then $M$ is Type~1 but not Type~3.
}
\end{example}
%
%
\begin{proposition}\label{3x3type}
Consider a differential field $\mathbb F$ of functions with field of constants
$\mathbb K$.
Let $M=[m_{i,j}]\in {\mathbb F}^{3,3}$ be $\mathbb K$-triangularizable and satisfy $M\, M'=M'\, M$. Then $M$ is of Type~\ref{1}.
\end{proposition}
\proof
Since $M$ is $\mathbb K$-triangularizable, we may assume that it is upper triangular already and consider different cases for the diagonal elements. If $M$ has three distinct diagonal elements, then it is $\mathbb K$-diagonalizable and the result follows by
Corollary~\ref{cor:diag_type1}.  If $M$ has exactly two distinct diagonal elements,
then it can be transformed to a direct sum of a $2\times 2$ and $1\times1$ matrix and hence the result follows by Proposition~\ref{2x2}. If all diagonal elements are equal, then, letting $E_{i,j}$ be the matrix that is zero except for the position $(i,j)$, where it is $1$, we have
$M=m_{1,1} I+ m_{1,3} E_{1,3}+ \widetilde M$, where $\widetilde M=m_{1,2} E_{1,2}+ m_{2,3} E_{2,3}$ also satisfies (\ref{c1}). Then it follows that $m_{1,2} m'_{2,3}= m'_{1,2} m_{2,3}$. If either
$m_{1,2}=0$ or $m_{2,3}=0$,  then we immediately have again Type~1, since $\widetilde M$ is a direct sum of a $2\times 2$ and a $1\times 1$ problem.
If both are nonzero, then $\widetilde M$ is nonderogatory and the result follows
by Theorem~\ref{th:nondero}. In fact, in this case $m_{1,2}=c m_{2,3}$ for some $c\in \mathbb K$ and therefore
\[
M=m_{1,1} I+ m_{1,3} E_{1,3}+ m_{2,3}\mat{ccc} 0 & c & 0 \\ 0 & 0 & 1 \\ 0 & 0 & 0 \rix,
\]
which is clearly of Type~\ref{1}.
\eproof
%

In the $4\times 4$ case, if the matrix is $\mathbb K$-triangularizable, then we either have at least two different eigenvalues, in which case we have reduced the problem again to the case of dimensions smaller than $4$, or there is only one eigenvalue,
and thus without loss of generality  $M$ is nilpotent. If $M$ is nonderogatory then we again have Type
$1$. If $M$ is derogatory then it is the direct sum of blocks of smaller dimension. If these dimensions are smaller than $3$, then we are again in the Type~1
case. So it remains to study the case of a block of size $3$ and a block of size $1$.
Since $M$ is nilpotent, the block of size $3$ is either Type~1 or Type~2. In both cases the complete matrix is also Type~1 or Type~2, respectively.

The following example shows that $\mathbb K$-triangularizability is not enough
to imply that the matrix is Type~1.
\begin{example}{\rm
Consider the $9\times 9$ block matrix
\[
\hat M=\mat{ccc} 0 & M_a & 0 \\ 0 & 0 & M_a \\ 0 & 0 & 0 \rix,
\]
where $M_a$ is the Type~2 matrix from Example~\ref{eva3x3}.
Then $\hat M$ is nilpotent upper  triangular and not of Type~1, 2, or 3, the latter two
facts due to its $\mathbb F$-rank being $2$.
}
\end{example}

 Already in the $5\times 5$
case, we can find examples that are none of the (proper) types.

\begin{example}{\rm
Consider $M=T^{-1} \diag(M_1,M_2) T$ with $T\in \mathbb K^{n,n}$,
$M_1\in \mathbb F^{3,3}$ of Type~2 (e.g., take $M_1=M_a$ as in Example~\ref{eva3x3})
and $M_2=\mat{cc} 0&1\\ 0& 0\rix$. Then clearly $M$ is
not of Type~1 and it is not of Type~2, since
it has an $\mathbb F$-rank larger than $1$. By definition it is not of Type~3 either.
Clearly examples of any size can be constructed by building direct sums
of smaller blocks.
}
\end{example}

Schur's letter states that for $n\geq 6$ there are other types. The following example demonstrates this.
\begin{example}{\rm Let $M_a$ be the Type~2 matrix in Example~\ref{eva3x3}
and form the block matrix
\[
A =\mat{cc} M_a & I \\ 0 & M_a\rix .
\]
 Direct computation shows $AA' = A'A$ but $A'A'' \not = A''A$.  Furthermore  $A^3 = 0$ and $A$ has $\mathbb F$-rank $3$. Thus $A$ is neither Type~1, Type~2  nor Type 3 (the last case need not be considered, since $A$ is nilpotent). We also note that $\rank(A'') = 6$. We now assume that $\mathbb K$ is algebraically closed  and we show that $A$ is not $\mathbb K$-similar to the direct sum of Type~1 or Type~2 matrices.

 To obtain a contradiction we assume that (after a $\mathbb K$-similarity)
$A = \diag(A_1,A_2)$ where $A_1$ is the direct sum of
Type~1 matrices (and hence Type~1) and $A_2$ is the direct sum of Type~2 matrices that are not Type~1. Since $A$ is not Type~1, $A_2$ cannot be the empty matrix. Since the minimum size of a Type~2 matrix that is not Type~1  is $3$ and its rank is $1$ it follows that  $A$ cannot be the sum of Type~2 matrices that are not Type~1. Hence the size of $A_1$ must be larger or equal to $1$ and, since $A_1$ is nilpotent, it follows that $\rank(A_1) < \size(A_1)$. Since
$A_1$ is $\mathbb K$-similar  to a strictly triangular matrix, it follows that
$\rank(A_1'')  < \size(A_1)$. Hence $\rank(A'') = \rank(A_1'') + \rank(A_2'') < 6$, a contradiction.
}
\end{example}

\begin{example}{\rm
If the matrix $M=\sum_{i=0}^r C_i t^i\in \mathbb F^{n,n}$ is a polynomial with coefficients $C_i\in \mathbb K^{n,n}$, then from (\ref{c1}) we obtain a specific set of conditions on
sums of commutators that have to be satisfied. For this we just compare coefficients
of powers of $t$ and  obtain a set of quadratic equations in the $C_i$, which has a
clear pattern. For example, in the case $r=2$, we obtain the three conditions $C_0C_1-C_1C_0=0$, $C_0C_2-C_2C_0=0$ and $C_1C_2-C_2C_1=0$, which shows that $M$ is of Type~1.
For $r=3$ we obtain the first nontrivial condition $3(C_0C_3-C_3C_0)+(C_1C_2-C_2C_1)=0$.

We have implemented
a Matlab routine for Newton's method to solve the set of quadratic matrix equations in the case $r=3$  and ran it for many different random starting coefficients
$C_i$ of different dimensions $n$. Whenever Newton's method converged (which it did in most of the
cases) it converged to a matrix of Type~1. Even in the neighborhood of a Type~2 matrix it converged to a Type~1 matrix. This suggests that the matrices of Type~1 are generic in the set of matrices satisfying~(\ref{c1}). A copy of the Matlab routine is available from the authors upon request.
}
\end{example}

\section{Conclusion}\label{conclusion}

We have presented a letter of Schur's that contains a major contribution
to the question when a matrix with elements that are functions in one variable
commutes with its derivative. Schur's letter precedes many partial results
on this question, which is still partially open. We have put Schur's result in perspective with later results and extended it in an algebraic context to matrices over a differential field. In particular, we have presented several results that characterize
Schur's matrices of Type~1. We have given examples of matrices that commute with their derivative which are of none of the Types~1,~2 or~3.We have shown that matrices of Type~1 may be triangularized over the constant field (which implies that their eigenvalues lie in the differential field) but we are left with an open problem already mentioned in Section \ref{sec:discussion}.
\begin{openproblem} {\rm Let $M$ be a matrix in a differential field $\mathbb F$, with an algebraically closed
field of constants, that satisfies $MM' = M'M$. Must the eigenvalues of $M$ be elements of the field $\mathbb F$?}
\end{openproblem}
For example, if $M$ is a polynomial matrix over the complex numbers must the eigenvalues be rational functions?
We have found no counterexample.

\section*{Acknowledgements}
We thank Carl de Boor for helpful comments on a previous draft of the paper and Olivier Lader for his careful reading of the paper and for his suggestions. We also thank an anonymous referee for pointing out the observation in Remark~\ref{refrem}.

\bibliographystyle{plain}

\end{document}